\documentclass{article}

\usepackage[english]{babel}

\usepackage[utf8]{inputenc}
\usepackage[T1]{fontenc}

\usepackage{amssymb}
\usepackage{stmaryrd}
\usepackage{amsmath,amsfonts}
\usepackage[tikz]{bclogo}
\usepackage[top=4.34cm]{geometry}

\usepackage{lmodern}
\usepackage{textcomp}
\usepackage{mathtools, bm}
\usepackage{amssymb, bm}
\usepackage[unq]{unique}
\usepackage{enumerate}
\usepackage{comment}
\usepackage[breaklinks]{hyperref}

\usepackage[numbers]{natbib}  

\usepackage[breaklinks]{hyperref}
\usepackage[numbers]{natbib}    

\usepackage{amsthm}

\usepackage{cleveref}

\newtheorem{theorem}{Theorem}[section]
\newtheorem{lemma}[theorem]{Lemma}

\newtheorem{conjecture}{Conjecture}[section]

\newtheorem{Claim}[theorem]{Claim}

\theoremstyle{definition}

\title{The odd chromatic number of a planar graph is at most $8$}

\author{Jan Petr\footnote{\href{mailto:jp895@cam.ac.uk}{jp895@cam.ac.uk}, Centre for Mathematical Sciences,
Wilberforce Road,
Cambridge CB3 0WA,
United Kingdom} \and Julien Portier\footnote{\href{mailto:jp899@cam.ac.uk}{jp899@cam.ac.uk}, Centre for Mathematical Sciences,
Wilberforce Road,
Cambridge CB3 0WA,
United Kingdom}}

\date{}

\begin{document}

\maketitle

\begin{abstract}

Petruševski and Škrekovski \cite{odd9} recently introduced the notion of an odd colouring of a graph: a proper vertex colouring of a graph $G$ is said to be \emph{odd} if for each non-isolated vertex $x \in V(G)$ there exists a colour $c$ appearing an odd number of times in $N(x)$.

Petruševski and Škrekovski proved that for any planar graph $G$ there is an odd colouring using at most $9$ colours and, together with Caro \cite{oddremarks}, showed that $8$ colours are enough for a significant family of planar graphs. We show that $8$ colours suffice for all planar graphs.
\end{abstract}

\section{Introduction}

Let $G$ be a graph. An \emph{odd colouring of $G$} is a proper colouring $\varphi$ of $G$ such that for every non-isolated vertex $x \in V(G)$ there is a colour $c$ satisfying that $|\varphi^{-1}(c) \cap N(x)|$ is odd. The \emph{odd chromatic number of $G$}, denoted $\chi_o(G)$, is the smallest number $k$ such that there exists an odd colouring using $k$ colours.

Odd colourings were recently introduced by Petruševski and Škrekovski \cite{odd9}. Using the discharging method, they showed that $\chi_o(G)\leq 9$ holds for all planar graphs $G$. Furthermore, they made the following conjecture.

\begin{conjecture}
Every planar graph $G$ has odd chromatic number at most 5.
\end{conjecture}

If true, this conjecture is clearly best possible as $\chi_o(C_5)=5$. Since the introduction of the notion of odd colouring in \cite{odd9}, two papers have appeared on arXiv containing partial results about narrowing down the gap between the bounds of $5$ and $9$ colours for planar graphs. Cranston \cite{oddsparse} focused on the odd chromatic number of sparse graphs, obtaining for instance that $\chi_o(G)\leq 6$ for planar graphs $G$ of girth at least $6$ and $\chi_o(G)\leq 5$ for planar graphs $G$ of girth at least $7$. Caro, Petruševski and Škrekovski \cite{oddremarks} studied various properties of the odd chromatic number, including the following important steps towards proving that $8$ colours suffice for an odd colouring of any planar graphs.

\begin{lemma}\label{8even}
If G is a connected planar graph of even order, then $\chi_{o}(G) \leq 8$.
\end{lemma}

\begin{lemma}\label{8odd}
If $G$ is a connected planar graph of odd order which has a vertex of degree $2$
or any odd degree, then $\chi_{o}(G) \leq 8$.
\end{lemma}

For the reader's convenience, we include the proofs by Caro, Petruševski and Škrekovski of those lemmata at the beginning of the next section. Their proofs use the following theorem by Aashtab, Akbari, Ghanbari and Shidani \cite{4trees}.

\begin{theorem}\label{4forest}
Let $G$ be a connected planar graph of even order. Then $V(G)$ partitions into at most $4$ sets such that each part induces an odd forest.
\end{theorem}

It is worth noting that the proof of \Cref{4forest} is based on the Four-Colour Theorem.

Since a minimal planar graph (with respect to the number of vertices) which does not admit an odd colouring with $8$ colours must be connected, \Cref{8even} and \Cref{8odd} obviously imply the following.

\begin{lemma}\label{ReduceCase}
Let $G$ be a minimal planar graph which does not admit an odd colouring with $8$ colours. Then $G$ has odd order and all degrees in $G$ are even and at least $4$.
\end{lemma}

The goal of this paper is to prove that $8$ colours are sufficient for an odd colouring of any planar graphs.

\begin{theorem}\label{MainThm}
For every planar graph $G$ we have $\chi_{o}(G) \leq 8$.
\end{theorem}

As usual, we call a vertex of degree $d$ a \emph{$d$-vertex}, and a vertex of degree at least $d$ a \emph{$d^{+}$-vertex}. A face is of size $k$ if the enumeration of the vertices on its border has length $k$, counting with multiplicity a vertex that appears multiple times. A face of size $k$ is called a \emph{$k$-face} and a face of size at least $k$ a \emph{$k^+$-face}.

\section{The proof}

We begin by giving the proofs of \Cref{8even} and \Cref{8odd} due to Caro, Petruševski and Škrekovski.

\begin{proof}[Proof of \Cref{8even}.]
By \Cref{4forest}, we partition $V(G)$ into at most $4$ sets such that each part induces an odd forest. Each odd forest $F$ obviously satisfies $\chi_0(F) = 2$. By using different colours for the different forests, we can find an odd colour $G$ with at most $8$ colours.
\end{proof}

Before moving onto the proof of \Cref{8odd}, we introduce a convenient notation that will be used many times in the rest of the paper: suppose we remove a vertex $v$ from a graph $G$ and possibly add some edges between neighbours of $v$ while keeping the resulting graph $G'$ planar. Assume $G'$ has an odd colouring using at most $8$ colours. To extend the colouring to $G$, we need to check two things: 1) that there is a colour that appears an odd number of times in $N(v)$, and 2) that there exists a colour for $v$ such that the colouring remains proper and the neighbours of $v$ still each have a colour that appears an odd number of times in their own neighbourhood. Note that for the second condition, every neighbour $w$ of $v$ prevents $v$ from using at most $2$ colours: the colour $c(w)$ (since we are looking for a proper colouring of $G$), and at most one other in the case where $w$ has exactly one colour that appears an odd number of times in $N(w)\setminus \{v\}$ in $c$. In the following, we will call the first of the colours \emph{forbidden by the proper colouring condition} and the second \emph{forbidden by the oddness condition}. 

\begin{proof}[Proof of \Cref{8odd}.]
If $v \in V(G)$ is a vertex of odd degree, we attach to it a new leaf $w$. The new graph has even order, and consequently admits an odd colouring by \Cref{8even}. Deleting $w$ gives an odd colouring of the original graph $G$ since $v$ is of odd degree and hence always has its oddness condition satisfied.

If $v \in V(G)$ is a vertex with $d(v)=2$, we remove $v$ and add an edge between the neighbours of $v$ (if they do not already have an edge between them). The new graph has even order, and consequently admits an odd colouring by \Cref{8even}. In the original graph $G$, $v$ has now $2$ colours forbidden by the proper colouring condition, and at most $2$ colours forbidden by the oddness condition. Hence there is a colour that we can give to $v$ to extend the colouring to an odd colouring of $G$.
\end{proof}

From now on, assume for contradiction that there is a planar graph which does not admit an odd colouring with at most $8$ colours, and let $G=(V,E)$ be a minimal such graph. By \Cref{ReduceCase}, all the degrees in $G$ are even and at least $4$.

We notice that $G$ cannot contain certain configurations of vertices:

\begin{Claim}\label{No4sWCommonNbr}
There do not exist two neighbouring vertices of degree $4$ in $G$, that share at least one common neighbour.
\end{Claim}

\begin{proof}
Suppose $u$ and $v$ are two neighbouring $4$-vertices with a common neighbour $w$. Consider the graph $G'$ obtained from $G$ by removing $u$ and $v$. Let $c$ be an odd colouring of $G'$ using at most $8$ colours.

\begin{figure}[htbp]\centering
    			\includegraphics[height=2cm]{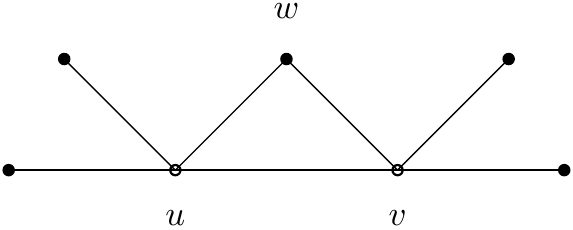}
    			\caption{An illustration of \Cref{No4sWCommonNbr}. All neighbours of $u$ and $v$ are drawn.}
\label{4and4}
	\end{figure}

First, suppose at least one of $u$ and $v$, say $u$, has another neighbour of colour $c(w)$. At most $7$ colours are forbidden for $v$: at most $3$ colours are used for $v$'s neighbours (we have not coloured $u$ yet), and at most $4$ colours are forbidden by the oddness condition. We then colour $v$ with any of the remaining colours. Now the neighbours of $u$ forbid at most $7$ colours, since the colour $c(w)$ is forbidden by the two neighbours having colour $c(w)$, hence there is a choice of a colour for $u$ that gives a valid odd colouring of $G$ with $8$ colours, which is a contradiction.

If neither $u$ nor $v$ has another neighbour of colour $c(w)$, then we colour $v$ in the same way as above. Now for $u$, the colour $c(w)$ is forbidden by two different sources: from $w$ by the proper colouring condition, and from $v$ by the oddness condition, leaving one colour available for $u$, which again gives a valid odd colouring of $G$ with $8$ colours, which is a contradiction.
\end{proof}

\begin{Claim}\label{No46neighbours}
There do not exist two neighbouring vertices $v,x$ in $G$ such that $v$ is a 4-vertex, $x$ a 6-vertex and $|N(v) \cap N(x)| \geq 2$.
\end{Claim}

\begin{proof}
Let $y_1$ and $y_2$ be two of the common neighbours of $v$ and $x$. Let $w$ be the fourth neighbour of $v$ and let $z_1,z_2,z_3$ be the other neighbours of $x$. See \Cref{4and6} for illustration.

\begin{figure}[htbp]\centering
    			\includegraphics[height=5cm]{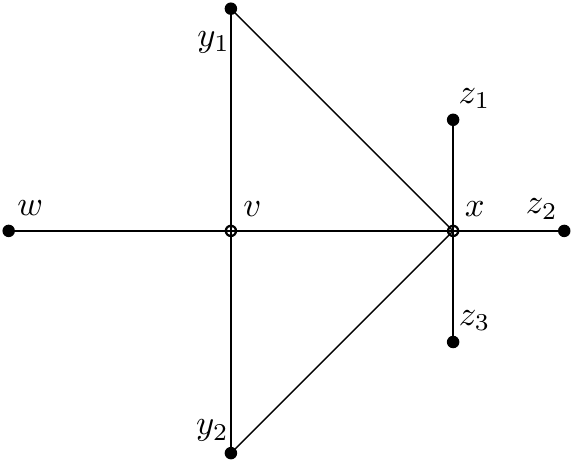}
    			\caption{An illustration of \Cref{No46neighbours}. All neighbours of $v$ and $x$ are drawn.}
\label{4and6}
	\end{figure}

Set $V'=V(G) \setminus v$ and $E'=E \cup \{wx\} \setminus \{vy_1, vx, vy_2, vw \}$, and let $c$ be a colouring of $G'=(V', E')$ using at most $8$ colours. For a vertex $u$, let $b(u)$ be the colour, if it exists, forbidden at $u$ by the oddness condition. If this colour does not exist, let $b(u)$ be an arbitrary colour. Note that the oddness condition is already satisfied for $v$ with $c(x)$, hence $c(y_1)$, $c(x)$, $c(y_2)$, $c(w)$, $b(y_1)$, $b(x)$, $b(y_2)$ and $b(w)$ must all be different, otherwise, there is a colour for $v$ that extends $c$ to an odd colouring on $G$.
To complete our proof of \Cref{No46neighbours} we will now find an odd colouring of $G$ agreeing with $c$ on $V \setminus \{v,x\}$, and so contradicting the definition of $G$.

If $b(y_1)$ and $b(y_2)$ both appear among $ b(z_1)$, $b(z_2)$ and $b(z_3)$, assign to $v$ the colour $c(x)$. Then $x$ has $4$ forbidden colours by the proper colouring condition: $c(y_1)$, $c(y_2)$, $c(u)$ and $c(x)$. But it has at most $3$ forbidden colours by the oddness condition: $b(z_1)$, $b(z_2)$ and $b(z_3)$. Thus there is a colour that can be assigned to $x$. It is clear that this colouring is proper. Moreover, since $c(x) \neq b(w)$, $w$ has a colour appearing an odd number of times in its neighbourhood. The oddness condition is also satisfied for all other vertices. Thus this colouring is an odd colouring.

If at least one of $b(y_1)$ and $b(y_2)$ does not appear among $b(z_1)$, $b(z_2)$ and $b(z_3)$, call this colour $b_1$, and the other one $b_2$. Then we can assign $b_1$ to $x$ and $b_2$ to $v$. It is straightforward to check that also in this case the resulting colouring is an odd colouring of the whole $G$, giving a contradiction.
\end{proof}

On top of the previous two claims we will use the discharging method. The general idea of the discharging method is to obtain certain local configurations in any planar drawing of the given graph. This is done by assigning initial charge to the vertices and faces so that the sum of all the charges is negative, and then distributing the charge according to so-called discharging rules. Those typically make many final charges non-negative. Since the sum of all charges does not change, there has to exist a vertex or a face whose final charge is negative. The discharging rules thus give information about the surroundings of this vertex or face.

In our case, we assign to each vertex $v$ initial charge $ch_0(v)=d(v)-6$ and to each face $F$ initial charge $ch_0(F)=2d(F)-6$ where by $d(v)$ we mean the degree of $v$ and by $d(F)$ we mean the number of edges on the boundary of $F$. Let $ch_1(v)$ and $ch_1(F)$ be the final charges of vertex $v$ and of face $F$. Applying Euler's formula, we have $\sum_{v} ch_1(v) + \sum_{F} ch_1(F) = \sum_{v} ch_0(v) + \sum_{F} ch_0(F)=-12$.

The discharging rules we apply are:

\begin{enumerate}[(R1)]
    \item Every $8^+$-vertex $u$ sends charge $1/2$ to every neighbouring $4$-vertex such that the next neighbour of $u$ in the counter-clockwise order is a $6^+$-vertex.
    \item Every $4^+$-face $F$ sends charges to its $4$-vertices according to the following rules:
    \begin{enumerate}[(i.)]
        \item If $F$ is a $4$-face with all vertices having degree $4$, then $F$ sends charge $1/2$ to each of its vertices.
        \item If $F$ is a $4$-face with three incident $4$-vertices and one incident $6^+$-vertex $u$, then $F$ sends charge $3/4$ to its $4$-vertices that are neighbours of $u$ and $1/2$ to the remaining $4$-vertex.
        \item If $F$ is a $5$-face with all vertices having degree $4$, then $F$ sends charge $3/4$ to all its $4$-vertices.
        \item Otherwise, $F$ sends charge $1$ to all its $4$-vertices.
    \end{enumerate}
\end{enumerate}

\begin{figure}[htbp]\centering
    			\includegraphics[height=4.5cm]{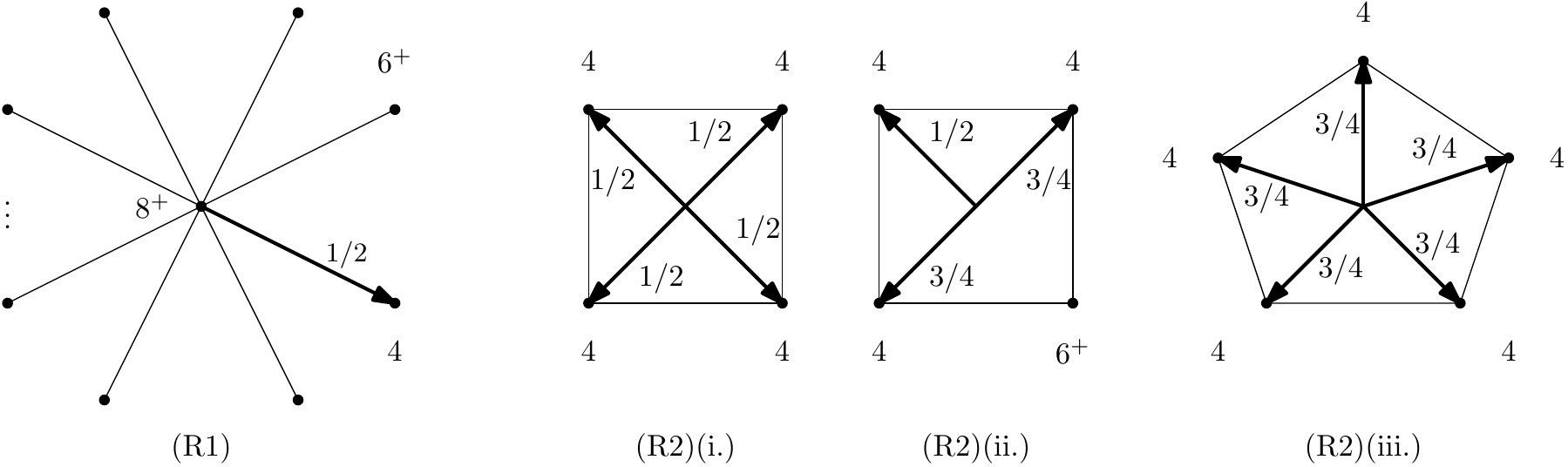}
    			\caption{The discharging rules. If a non-triangular face does not look like one of the three depicted above, it sends charge $1$ to each of its incident $4$-vertices, in accordance with (R2)(iv.).}
\label{discharge}
	\end{figure}
	
As in the definition of a $k$-face we gave earlier, if a vertex appears multiple times on the border of the face, then the face sends it charges with multiplicity.

The discharging rules were chosen so that only a $4$-vertex can have a negative final charge. Before verifying this, note also that each $4$-vertex receives charge at least $1/2$ from each non-triangular face it is incident to.

\begin{Claim} \label{vertices}
After discharging, every $6^+$-vertex carries a non-negative charge.
\end{Claim}
\begin{proof}
Only (R1) affects the charges of $6^+$-vertices. The charge of $6$-vertices remains $0$ after discharging. Let $x$ be an $8^+$-vertex. It has at most $\frac{d(x)}{2}$ neighbouring $4$-vertices such that the previous neighbour in the counter-clockwise order is a $6^+$-vertex. Therefore, we have: $$ch_1(x) \geq (d(x)-6)-\frac{1}{2} \cdot \frac{1}{2} \cdot d(x) \geq 0.$$ \end{proof}

\begin{Claim} \label{faces}
After discharging, each face carries a non-negative charge.
\end{Claim}
\begin{proof}
Only (R2) can change the charge of a face. The charge of a $3$-face remains $0$. Let $F$ be a $6^+$-face. Then, by (R2)(iv.), $ch_1(F) \geq 2d(F)-6-d(F) \geq 0$.

If $F$ is a $5$-face, there are two possibilities. If all its incident vertices are $4$-vertices, then the charge $F$ sends is $5 \cdot 3/4$ as by (R2)(iii.), and thus $ch_1(F)=2 \cdot 5-6-15/4>0$. On the other hand, if there is at least one $6^+$-vertex incident with $F$, then $F$ sends charge $1$ to at most $4$ vertices. Then, $ch_1(F) \geq 2 \cdot 5-6-4 \cdot 1=0$.

Finally, let $F$ be a $4$-face. If there are at most two $4$-vertices incident with it, $F$ sends a charge at most $2$, and thus $ch_1(F)\geq 2 \cdot 4-6-2=0$. If there are exactly three $4$-vertices incident with $F$, then $F$ sends charge $2$ in accordance with (R2)(ii.). The face $F$ also sends charge $2$ if all four vertices incident with it are $4$-vertices in accordance with (R2)(i.). Thus, in every case, $ch_1(F)\geq 0$. 
\end{proof}

We are now able to deduce our main theorem from the claims above.

\begin{proof}[Proof of \Cref{MainThm}.]
Assume for contradiction that there exists a planar graph with odd chromatic number at least $9$. Let $G=(V,E)$ be such a graph that has minimal number of vertices. In view of \Cref{ReduceCase}, all degrees in $G$ are even and at least $4$. We fix a planar drawing of $G$ and apply the discharging rules (R1) and (R2). As $\sum_{v} ch_1(v) + \sum_{F} ch_1(F)=-12$, by \Cref{vertices} and \Cref{faces} there exists a negatively charged $4$-vertex $v$. Note that $ch_0(v) = -2$.

Let the neighbours of $v$ in clockwise order be $x_1,x_2,x_3,x_4$. The four faces around $v$ containing $\{x_1,x_2\}$, $\{x_2,x_3\}$, $\{x_3,x_4\}$ and $\{x_4,x_1\}$ are all different. Indeed, if any two of them had been the same, the face would have been a $6^+$-face, sending charge at least $2$ to $v$.

Since $ch_1(v)<0$ and $ch_0(v)=-2$, by (R2) at least one of these faces is a $3$-face, call it $F$. We can assume without loss of generality that $x_1x_2 \in E$. By \Cref{No4sWCommonNbr}, both $x_1,x_2$ are $6^+$-vertices.

By (R2) again, there is at least one other $3$-face around $v$ different than the face $(x_1x_2v)$, as otherwise the faces around $v$ containing $\{x_1,x_4\}$ and $\{x_2,x_3\}$ would each send charge at least $3/4$ to $v$ in addition to the charge at least $1/2$ sent to $v$ by the remaining non-triangular face.

We claim that one of $x_1x_4,x_2x_3$ is an edge. Indeed, if the above-mentioned $3$-face is formed by $v,x_3,x_4$, then $x_3$ and $x_4$ are $6^+$-vertices by \Cref{No46neighbours} which in turn means that one of the remaining two faces has to be a $3$-face too, as otherwise in view of (R2) each of them would send charge $1$ to $v$. Without loss of generality let $x_2x_3 \in E$. By \Cref{No4sWCommonNbr}, $x_3$ is a $6^+$-vertex.

In view of \Cref{No46neighbours}, $x_2$ is an $8^+$-vertex and thus sends charge $1/2$ to $v$ in accordance with (R1). By (R2), at least one of $x_3x_4,x_4x_1$ is an edge, as otherwise the two non-triangular faces around $v$ each send charge $3/4$ to it. In either case, $x_4$ is a $6^+$-vertex by \Cref{No4sWCommonNbr}. Applying \Cref{No46neighbours} again, either $x_3$ is an $8^+$-vertex (in case $x_3x_4 \in E$) or $x_1$ is an $8^+$-vertex (if $x_4x_1 \in E$).

Now, $v$ receives charge at least $1$ due to (R1). Therefore, the remaining face around $v$ is a $3$-face, as otherwise it would send charge $1$ in accordance with (R2). This means that $x_1x_2,x_2x_3,x_3x_4,x_4x_1$ are all edges and the respective faces around $v$ are all $3$-faces. Applying \Cref{No4sWCommonNbr} and \Cref{No46neighbours}, all of $x_1,x_2,x_3,x_4$ are $8^+$-vertices. But by (R1), $v$ receives charge $1/2$ from all of them, contradicting $ch_1(v)<0$.

Therefore, there is no planar graph with odd chromatic number greater than $8$.
\end{proof}

\section*{Acknowledgement}

The authors would like to thank their PhD supervisor Professor Béla Bollobás for his valuable comments.

\bibliographystyle{abbrvnat}  
\renewcommand{\bibname}{bib}
\bibliography{bib}

\end{document}